\documentclass{amsart}
\usepackage{amssymb}
\usepackage{latexsym}
\usepackage{amsmath}
\usepackage{euscript}
\usepackage{graphics}
\usepackage[all]{xy}
\usepackage[centering]{geometry}

\usepackage{MnSymbol} 

\newcommand{\be}{\begin{equation}}
\newcommand{\ee}{\end{equation}}
\newcommand{\ba}{\begin{eqnarray}}
\newcommand{\ea}{\end{eqnarray}}
\newcommand{\baa}{\begin{eqnarray*}}
\newcommand{\eaa}{\end{eqnarray*}}
\newcommand{\bb}{\begin{thebibliography}}
\newcommand{\eb}{\end{thebibliography}}

\newcommand{\bi}[1]{\bibitem{#1}}
\newcommand{\lab}[1]{\label{#1}}
\newcommand{\re}[1]{(\ref{#1})}

\newcounter{my}
\newcommand{\he}%
   {\stepcounter{equation}\setcounter{my}%
   {\value{equation}}\setcounter{equation}0%
   }%
\newcommand{\she}%
   {\setcounter{equation}{\value{my}}%
    }%

\renewcommand\t{\tilde}

\newcommand\ve{\varepsilon}

\newtheorem{lem}{Lemma}

\theoremstyle{definition}

\numberwithin{equation}{section}

\begin{document}

\title[symmetric hypergeometric polynomials]{Symmetric abstract hypergeometric polynomials}


\author{Satoshi Tsujimoto}
\author{Luc Vinet}
\author{Guo-Fu Yu} 
\author{Alexei Zhedanov}

\address{Department of Applied Mathematics and Physics, Graduate School of Informatics, Kyoto University, Yoshida-Honmachi, Kyoto, Japan 606--8501}

\address{Centre de recherches math\'ematiques
Universite de Montr\'eal, P.O. Box 6128, Centre-ville Station,
Montr\'eal (Qu\'ebec), H3C 3J7}

\address{Department of Mathematics, Shanghai Jiao Tong University
Shanghai 200240, China}

\address{Department of Mathematics, School of Information, Renmin University of China, Beijing 100872,
CHINA}

\begin{abstract}
Consider an abstract operator $L$ which acts on monomials $x^n$  according to
$L x^n= \lambda_n x^n + \nu_n x^{n-2}$ for $\lambda_n$ and $\nu_n$ some coefficients. Let $P_n(x)$ be eigenpolynomials of degree $n$ of $L$: $L P_n(x) = \lambda_n P_n(x)$. A classification of all the cases for which the polynomials $P_n(x)$ are orthogonal is provided. A general derivation of the algebras explaining the bispectrality of the polynomials is given. The resulting algebras prove to be  central extensions of the Askey-Wilson algebra and its degenerate cases.
\end{abstract}

\keywords{}


\maketitle

\section{Introduction}
There is a much interest in obtaining abstractly a classification of orthogonal polynomials  (OPs) that possess certain general properties. This indeed allows to characterize many families of OPs simultaneously and to shed light on the structural origin of features that they share. We here pursue in this respect a program undertaken recently \cite{Zhe_hyp}, \cite{VZ_Newton}, \cite{Z_umbral} that aims to determine the orthogonal polynomials that are eigenfunctions of an operator $L$ that has a prescribed action in monomial or Newtonian bases.

The simplest case, discussed in \cite{Zhe_hyp}, is one where $L$ acts according to $L x^n = \lambda_n x^n + \rho_n x^{n-1}$ on monomials ($\mu_n$ and $\rho_n$ are coefficients). We here examine the situation when the operator $L$ rather acts as follows
\be
L x^n = \lambda_n x^n +   \nu_n x^{n-2}, \quad n=0,1,2,\dots \lab{main_L_x_n} \ee
with $\lambda_n$ and $\nu_n$ sets of real factors. All eigenpolynomials of $L$ that form orthogonal ensembles will be determined and found to be the generalized Gegenbauer and Hermite polynomials and their $q$-analogs. Expectedly, these are symmetric polynomials.

The algebras encoding the bispectrality of these families shall also be obtained in an abstract fashion and seen to involve the reflection operator as a central element. 

The outline of the paper is as follows. Some background is provided in Section 2 where the essential elements  of the analysis presented in \cite{Zhe_hyp} will be reviewed. These results will prove useful to obtain in Section 3 the classification of the (necessarily) symmetric eigenpolynomials of the operator $L$ defined by eq. \re{main_L_x_n}. Section 4 is dedicated to the differential - difference realizations of the operators L that emerge in the classification and Section 5 focuses on the construction of the nonlinear algebras associated to these polynomials. The paper ends with concluding remarks.

\section{Background}
\setcounter{equation}{0} 
Let $L$ be an operator which transforms any polynomial of degree $n=0,1,2,\dots$ to a polynomial of the same degree $n$. This operator preserves the $(n+1)$-dimensional linear space of polynomials with degrees $\le n$ and hence the eigenvalue problem for $L$ on that space is well-posed.  

To specify the operator $L$, we assume that its action on monomials is the following
\be
L x^n = \lambda_n x^n + \nu_n^{(1)} x^{n-1} + \nu_n^{(2)} x^{n-2} + \dots + \nu_n^{(J)} x^{n-J}, \lab{L_x_n} \ee
where $J=1,2,\dots$ is a fixed positive integer and where $\lambda_n, \nu_n^{(i)}$ are some coefficients. The eigenvalue problem will then read as follows in terms of monic polynomials $P_n(x)= x^n + O(x^{n-1})$: 
\be
L P_n(x) = \lambda_n P_n(x), \quad n=0,1,2,\dots \lab{eig_LP} \ee 
It is assumed  in what follows that the eigenvalues $\lambda_n$ are nondegenerate, i.e.
\be
\lambda_n \ne \lambda_m, \quad \mbox{if} \quad n \ne m . \lab{ndeg_lambda} \ee
It can easily be checked that under condition \re{ndeg_lambda}, the eigenpolynomial $P_n(x)$ defined by \re{eig_LP} is unique. Indeed, if one assumes that there exists another monic polynomial $\t P_m(x), \: m\le n$ with the same eigenvalue $\lambda_n$  it follows that $m=n$. Hence there exists two monic polynomials $P_n(x)$ and $\t P_n(x)$ with the same eigenvalue $\lambda_n$; their difference $\Delta P_n(x) = P_n(x) - \tilde P_n(x)$ is a polynomial of a smaller degree, say $j <n$ which again satisfy the eigenvalue equation $L \Delta P_n(x) = \lambda_n \Delta P_n(x)$. This is however impossible as it would contradict the condition \re{ndeg_lambda} thus showing that $P_n$ is unique.

The main problem is to determine when the polynomials $P_n(x), \: n=0,1,2,\dots$ form an orthogonal set. This is so if there exists a nondegenerate linear functional $\sigma$ acting on the space of polynomials such that
\be
\langle \sigma, P_n(x) P_m(x) \rangle = 0, \quad n \ne m \lab{ort_P} \ee
where $\langle \sigma, \pi(x) \rangle$ stands for the evaluation of the functional $\sigma$ on the polynomial $\pi(x)$. The functional $\sigma$ is completely determined by the moments
\be
\langle \sigma, x^n \rangle = c_n, \quad n=0,1,2,\dots \lab{mom_def} \ee
and is said nondegenerate if $\Delta_n \ne 0, \;
n=0,1,2,\dots$, where $\Delta_n = |c_{i+k}|_{i,k=0}^n$ are the Hankel determinants constructed from the moments. The conditions \re{mom_def} are equivalent to the conditions 
\be
\langle \sigma, P_n(x) \pi(x) \rangle = 0, \lab{ort_P2} \ee
where $\pi(x)$ is any polynomial of degree lesser than $n$. Equivalently, the polynomials $P_n(x)$ are orthogonal iff they satisfy the three-term recurrence relation \cite{Chi}
\be
P_{n+1}(x) + b_n P_n(x) + u_{n} P_{n-1}(x) = xP_n(x), \quad n=1,2,\dots \lab{rec_P} \ee 
with some  coefficients $b_n, u_n$. The linear functional $\sigma$ is nondegenerate iff $u_n \ne 0, \; n=1,2,\dots$. When the coefficients $b_n$ are real and $u_n>0, \; n=1,2,\dots$, the polynomials are orthogonal with respect to a positive measure $d \mu(x)$ on the real axis \cite{Ismail}:
\be
\int_{a}^b P_n(x) P_m(x) d \mu(x) = h_n \delta_{nm}, \lab{ort_mu} \ee
where $h_0=1$ and $h_n=u_1 u_2 \dots u_n, \: n=1,2,\dots$. The integration limits $a$, $b$ in \re{ort_mu} may be either finite or infinite.

It is shown in \cite{Zhe_hyp}  that the necessary and sufficient conditions for the polynomials $P_n(x)$ to be orthogonal are that 
\be
\langle \sigma, g(x) L f(x) \rangle
=  \langle \sigma, f(x) L g(x) \rangle \lab{gLf} \ee
for any pair of polynomials $f(x), g(x)$. 
This means that $L$ must be symmetric on the space of polynomials with respect to the functional $\sigma$. See also \cite{Duran} where the same condition is derived when $L$ is a higher-order difference operator. In fact, \re{gLf} implies the orthogonality of the eigenpolynomials   for {\it any} linear operator $L$ preserving polynomiality and such that $\deg(L p(x)) \le \deg(p(x))$ for any polynomial $p(x)$. Choosing $f(x) = x^n$ and $g(x)=x^m$ and assuming that the operator $L$ acts as in \re{L_x_n},  condition \re{gLf} is seen to amount to 
\be
(\lambda_n - \lambda_m)c_{n+m} + (\nu_n^{(1)} - \nu_m^{(1)})c_{n+m-1} + \dots + (\nu_n^{(J)} - \nu_m^{(J)})c_{n+m-J} =0, \lab{gen_NS} \ee
where $n,m$ can take all admissible values (i.e. $n+m \ge J$).

When $J=1$, \re{gen_NS} becomes
\be
(\lambda_n - \lambda_m)c_{n+m} + (\nu_n^{(1)} - \nu_m^{(1)})c_{n+m-1}=0. \lab{NS_J=1} \ee
From the no-degeneracy condition $\lambda_n \ne \lambda_m$ for all $n \ne m$, it follows \cite{Zhe_hyp} that  
\be
\frac{\lambda_{n+1}-\lambda_k}{\nu_{n+1}-\mu_k}=\frac{\lambda_{n}-\lambda_{k+1}}{\nu_{n}-\mu_{k+1}}
\lab{cross} \ee
should be valid for all pairs of integers $k,n=0,1,2,\dots$.

As shown in \cite{Zhe_hyp}, all nondegenerate solutions of \re{cross} are obtained as solutions of the linear difference equations: 
\be
\lambda_{n+1}+\lambda_{n-1} -\Omega  \lambda_n +B_1 =0, \quad
\nu^{(1)}_{n+1}+\nu^{(1)}_{n-1} -\Omega  \nu^{(1)}_n +B_2 =0, \lab{dfr_main} \ee
where $\Omega, B_1, B_2$ are arbitrary parameters. The type of solution depends on the choice of the parameter $\Omega$.

When $\Omega=q+q^{-1}>2$ with $q$ a real parameter, the general solution is
\be
\lambda_n = \alpha_1 q^n + \alpha_2 q^{-n} + \alpha_0, \quad \nu_n^{(1)} = \beta_1 q^n + \beta_2 q^{-n} + \beta_0 \lab{ln_0_q} \ee
where $\alpha_1, \alpha_2, \beta_1, \beta_2$ are arbitrary parameters and $\alpha_0, \beta_0$ are related to $B_1, B_2$ as follows
\be
B_1 =\alpha_0 q^{-1}(1-q)^2, \quad B_2 =\beta_0 q^{-1}(1-q)^2. \lab{alb0_B_q} \ee
This case corresponds to the little $q$-Jacobi polynomials and their ``specializations'': $q$-Laguerre, little $q$-Laguerre, alternative $q$-Charlier and Stieltjes-Wigert polynomials.

When $\Omega=2$ the general solution is
\be
\lambda_n = \alpha_2 n^2 + \alpha_1 n + \alpha_0 , \quad \nu_n^{(1)} =\beta_2 n^2 + \beta_1 n + \beta_0 \lab{ln_0_1} \ee
where $\alpha_0, \alpha_1, \beta_1, \beta_0$ are arbitrary constants and 
\be
B_1= -2 \alpha_2, \; B_2 = -2 \beta_2. \lab{alb0B_1} \ee
This case corresponds to the Jacobi and Laguerre polynomials.

Finally, when $\Omega=-2$,
\be
\lambda_n = \alpha_1 (-1)^n + \alpha_2 (-1)^n n + \alpha_0, \quad \nu_n^{(1)}=\beta_1 (-1)^n + \beta_2 (-1)^n n + \beta_0 \lab{ln_0_q=-1} \ee
with $\alpha_1, \alpha_2$ arbitrary and
\be
B_1 =-4 \alpha_0, \; B_2 = -4 \beta_0. \lab{alb0B_q=-1} \ee
This case corresponds to the little -1 Jacobi polynomials.

The case $\Omega<-2$  is basically obtained from the one when $\Omega > 2$ under the substitution $q \to -q$.
The remaining case $|\Omega|<2$ requires that q belongs to the unit circle $|q|=1$; the corresponding orthogonal polynomials are not positive definite \cite{Zhe_hyp} and will not be considered here.
It should be stressed that the parameters $\alpha_i$ and $\beta_i$ in the above formulas can be arbitrary.  However,  from the condition that $L(1)=\lambda_0$ is a constant it follows that necessarily $\nu_0^{(1)}=0$. This leads to obvious restrictions on the parameters $\alpha_i, \beta_i$.

\section{The classification problem and admissibility conditions}
\setcounter{equation}{0}
We have reviewed above the results of the classification of all the orthogonal eigenpolynomials of $L$ as given in \re{L_x_n} when $J=1$. In order of complexity, the next classification problem corresponds to $J=2$. In this case, the necessary and sufficient conditions  \re{gen_NS} for the polynomials $P_n(x)$ to be orthogonal specialize to 
\be
(\lambda_n - \lambda_m)c_{n+m} + (\nu_n^{(1)} - \nu_m^{(1)})c_{n+m-1} + (\nu_n^{(2)} - \nu_m^{(2)})c_{n+m-2} =0.
\lab{NS_J=2} \ee
This system of equations proves rather complicated and we shall restrict ourselves to the simpler case that result when $\nu_n^{(1)}=0$. This means that we will study the operator $L$ which acts as follows
\be
L x^n = \lambda_n x^n +   \nu_n x^{n-2} \lab{main_L_x} \ee
with two unknown coefficients $\lambda_n, \nu_n$ (we omit the superscript in $\nu_n$). We assume that $\lambda_n$ satisfy the nondegenerate condition \re{ndeg_lambda} and moreover that
\be
\nu_0 =\nu_1=0, \quad \nu_n \ne 0, \; n=2,3,4, \dots \lab{nu_0} \ee
The necessary and sufficient conditions for orthogonality will now read:
\be
 (\lambda_n - \lambda_m)c_{n+m} +  (\nu_n - \nu_m)c_{n+m-2} =0.
\lab{NS_main} \ee
The classification problem that we will tackle now amounts to finding all the orthogonal polynomials that are eigenfunctions of an operator of the form specified by \re{main_L_x}.

We first prove the following simple lemma.
\begin{lem}
The monic eigenpolynomials $P_n(x)=x^n + O(x^{n-1})$ are uniquely determined by the equation $L P_n(x) =\lambda_n P_n(x)$ if the eigenvalues $\lambda_n$ are nondegenerate. Moreover, these eigenpolynomials are symmetric, i.e. they satisfy the parity condition
\be
P_n(-x) = (-1)^n P_n(x). \lab{par_P} \ee
\end{lem}
We already proved that the polynomial $P_n(x)$ is unique if the eigenvalues $\lambda_n$ are nondegenerate. Let us introduce the parity operator $R$ that acts on any function $f(x)$ as
\be
Rf(x) = f(-x). \lab{R_def} \ee
For all monomials it is clear that 
\be
LR x^n =RL x^n = (-1)^n \left( \lambda_n x^n +   \nu_n x^{n-2} \right). \lab{com_mon} \ee
Hence the operator $R$ commutes with $L$
\be
LR=RL \lab{com_LR} \ee
on the space of all polynomials.

This means that the operator $R$ leaves any eigenpolynomial $P_n(x)$ in the same eigenspace. Because of the uniqueness of the polynomial $P_n(x)$, it follows that the various eigensubspaces are one-dimensional and hence 
\be
R P_n(x) = \pm P_n(x). \lab{RP} \ee
In turn, condition \re{RP} means that the polynomial $P_n(x)$ should be an even or an odd function. We hence arrive at  \re{par_P}.

This condition implies that every even polynomial $P_{2n}(x)$ contains only the even monomials $x^{2k}$ while every odd polynomial $P_{2n+1}(x)$ contains only the odd monomials $x^{2k+1}$:
\be
P_{2n}(x) = U_n(x^2), \quad P_{2n+1}(x) = x V_n(x^2), \lab{P_UV} \ee
where $U_n(x),V_n(x)$ are monic polynomials of degree $n$. Such orthogonal polynomials are called symmetric  \cite{Chi}.

Assume now that the polynomials $P_n(x)$ are orthogonal with corresponding moments $c_n$. Condition \re{par_P} means that all odd moments are zero $c_{2n+1}=0, \: n=0,1,\dots$. This is equivalent to the statement that in the recurrence relation \re{rec_P}, the diagonal recurrence coefficients $b_n$ vanish \cite{Chi} and thus
the symmetric orthogonal polynomials satisfy the recurrence relation
\be
P_{n+1}(x) + u_n P_{n-1}(x) = xP_n(x). \lab{rec_sym} \ee
Assume that the polynomials $P_n(x)$ have the following expansion 
\be
P_n(x) = \sum_{k=0}^n A_{nk} x^k. \lab{P_exp} \ee 
The coefficients $A_{nk}$ have the property that $A_{nk}=0$ when the parities of $n$ and $k$ are different. This follows from \re{P_UV}. Moreover, $A_{nn}=1$ because the polynomials $P_n(x)$ are monic. 

Substituting \re{P_exp} into the recurrence relation \re{rec_sym}, we find on the one hand, the system of conditions
\be
A_{n+1,k} + u_n A_{n-1,k} = A_{n,k-1} , \quad k=0,1,\dots,n-1 .\lab{A_u_cond} \ee
For $k=n-1$, we obtain from \re{A_u_cond} the expression for $u_n$:
\be
u_n= A_{n+1,n-1}-A_{n,n-2}. \lab{u_A} \ee
On the other hand, we have from the eigenvalue equation $L P_n(x) = \lambda_n P_n(x)$ that 
\be
A_{nk}(\lambda_{n} -\lambda_k) + \nu_k A_{n,k+2} =0, \quad k=n-2, n-4, n-6, \dots \lab{rec_A} \ee  
This relation allows to write explicitly the coefficients $A_{nk}$ in terms of $\lambda_n, \nu_n$:
\be
A_{n,n-2} = -\frac{\nu_n}{\lambda_{n}-\lambda_{n-2}}, \; A_{n,n-4} = \frac{\nu_n \nu_{n-2}}{(\lambda_{n}-\lambda_{n-2})(\lambda_n-\lambda_{n-4})},  \dots \lab{A_expl} \ee
In particular, comparing \re{A_expl} with \re{u_A},  we arrive at the following expression for the recurrence coefficient $u_n$  in terms of $\lambda_n, \nu_n$:
\be
u_n = \frac{\nu_n}{\lambda_{n}-\lambda_{n-2}} -\frac{\nu_{n+1}}{\lambda_{n+1}-\lambda_{n-1}}. \lab{u_expl} \ee
Moreover, an elementary analysis shows (see, e.g. \cite{Chi}) that the polynomials $U_n(x)$ and $V_n(x)$  are again orthogonal  having the moments
\be
c^{(U)}_n=c_{2n}, \quad c^{(V)}_n=c^{(U)}_{n+1}=c_{2n+2}. \quad n=0,1,2,\dots \lab{UV_mom} \ee
It also follows  that the polynomials $V_n(x)$ are the kernel polynomials of the polynomials $U_n(x)$ \cite{Chi}:
\be
V_n(x)= \frac{U_{n+1}(x) - A_n U_n(x)}{x}, \quad A_n= \frac{U_{n+1}(0)}{U_n(0)}. \lab{V_ker} \ee
Let us introduce the coefficients
\be
\lambda_n^{(U)} = \lambda_{2n}, \;  \lambda_n^{(V)} = \lambda_{2n+1}, \; \nu_n^{(U)} = \nu_{2n}, \; \nu_n^{(V)} = \nu_{2n+1}, \quad n=0,1,2,\dots \lab{lanu_UV} \ee
It follows from \re{NS_main}  that
\be
(\lambda_n^{(U)}-\lambda_m^{(U)})c^{(U)}_{n+m} + (\nu_n^{(U)}-\nu_m^{(U)})c^{(U)}_{n+m-1}=0 \lab{U_cond} \ee
and 
\be
(\lambda_n^{(V)}-\lambda_m^{(V)})c^{(V)}_{n+m} + (\nu_n^{(V)}-\nu_m^{(V)})c^{(V)}_{n+m-1}=0. \lab{V_cond} \ee
Comparing conditions \re{U_cond}, \re{V_cond} with \re{NS_J=1}, we see that both polynomials $U_n(x)$ and $V_n(x)$ are abstract hypergeometric orthogonal polynomials of the $J=1$ type and must hence belong to one of the three classes described in Section 2.      This yields explicit expressions for the moments $c^{(U)}_n, c^{(V)}_n$ and hence for the moments $c_n$. One should not loose sight however that according to \re{UV_mom},  $c_n^{U}$ and $c_n^{V}$ collectively provide the moments of the polynomials $P_n(x)$ and that in view of \re{NS_main}, we have the following additional restrictions:
\be
\frac{\lambda_n^{(V)}-\lambda_m^{(V)}}{\nu_n^{(V)}-\nu_m^{(V)}} = \frac{\lambda_{n+1}^{(U)}-\lambda_m^{(U)}}{\nu_{n+1}^{(U)}-\nu_m^{(U)}}. \lab{add_res} \ee
It is not difficult to determine the solutions to these conditions. The results can be summarized as follows:

(i) Both polynomials $U_n(x)$ and $V_n(x)$ belong to the little $q$-Jacobi class or to its special and degenerate cases (say, $q$-Laguerre polynomials); generically, for little $q$-Jacobi polynomials we choose:
\ba
&&\lambda_n^{(U)}=(1-q^{2n})(q^{-2n} +a)+\ve_0, \; \lambda_n^{(V)} = \rho (q^{2n}-1)(q^{-2n} +aq^2)+\ve_1, \nonumber \\
&& \nu_n^{(U)}=(q^{2n}-1)(q^{-2n} +b), \: \nu_n^{(V)}=\rho(1-q^{2n})(q^{-2n} +bq^2)  \lab{la_nu_jq} \ea
where $\rho, \ve_0, \ve_1$ are arbitrary real parameters ($\rho$ is assumed to be nonzero). Note that the expressions for $\lambda^{(V)}_n, \nu^{(V)}_n$ are slightly different from those for $\lambda^{(U)}_n, \nu^{(U)}_n$; this is due to condition \re{add_res}. 

The meaning of the free parameters $\rho, \ve_0, \ve_1$ is the following. Assume that $L$ is an admissible operator, i.e. that the coefficients $\lambda_n$ and $\nu_n$ correspond to symmetric orthogonal polynomials $P_n(x)$ which are eigenfunctions of the operator $L$. We can introduce a more general operator
\be
\t L = \xi_1 L + \xi_2 L R + \eta_1 \mathcal{I} + \eta_2 R, \lab{tL} \ee
where $\mathcal{I}$ is the identity operator. The reflection operator $R$ commutes with $L$. Hence the new operator $\t L$ will again be an admissible operator for the same symmetric orthogonal polynomials $P_n(x)$ but with different eigenvalue:
\be
\t L P_n(x) = \t \lambda_n P_n(x) \lab{tLP} \ee
where 
\be
\t \lambda_n = (\xi_1 + (-1)^n \xi_2)\lambda_n + \eta_1 + \eta_2 (-1)^n. \lab{t_lambda} \ee
Among these 4 parameters, $\xi_1$ can be considered as a common factor and we can always put $\xi_1=1$. The other three constants $\xi_2,\eta_1,\eta_2$, correspond to the three  parameters $\rho, \ve_0, \ve_1$ in \re{la_nu_jq}. In what follows we will use the freedom in the  choice of these parameters in order to obtain the most simple presentations of the algebra associated to the various families of symmetric polynomials (see Section 5).

From \re{la_nu_jq} and \re{u_expl}, the explicit expression for the recurrence coefficient $u_n$ are found to be:
\be
u_{2n} = {\frac {{q}^{2\,n} \left( 1-{q}^{2\,n} \right)  \left( a{q}^{2\,n-2}-b
 \right) }{ \left( 1+a{q}^{4\,n} \right)  \left( 1+a{q}^{4\,n-2}
 \right) }}, \quad u_{2n+1} = {\frac {{q}^{2\,n} \left( b{q}^{2\,n+2}+1 \right)  \left( a{q}^{2\,n}
+1 \right) }{ \left( 1+a{q}^{4\,n} \right)  \left( 1+a{q}^{4\,n+2}
 \right) }}.
 \lab{u_qj} \ee
As expected, these recurrence coefficients correspond to those of the generalized little $q$-Jacobi polynomials \cite{GK}.

Instead of classifying all degenerate cases of the generalized little $q$-Jacobi polynomials we consider only the one with
\ba
&&\lambda_n^{(U)}=(1-q^n) +\ve_0, \; \lambda_n^{(V)} = \rho (1-q^n)+\ve_1,  \nonumber \\
&& \nu_n^{(U)}=(1-q^n)(q^{-n} -b), \; \nu_n^{(V)}=\rho (1-q^n)(q^{-n} -bq) \lab{la_nu_lq} \ea
As in the previous situation, we have for the recurrence coefficients the expression
\be
u_{2n} = b{q}^{2-2\,n} \left( {q}^{2\,n}-1 \right) , \quad u_{2n+1} = {q}^{-4\,n-2} \left( {q}^{2}+1-{q}^{2\,n+2} \right)  \left( b{q}^{2\,n
+2}+1 \right). \lab{u_qlag} \ee
These recurrence coefficients correspond to those of the generalized $q$-Laguerre polynomials \cite{GK}.

\medskip

(ii) Both polynomials $U_n(x)$ and $V_n(x)$ are either Jacobi or Laguerre polynomials. Namely, for the Jacobi polynomials
\ba
&&\lambda_n^{(U)}=-2n(2n+a) + \ve_0, \qquad \lambda_n^{(V)} = 2\rho n(2n+a+2)+\ve_1, \nonumber \\
&& \nu_n^{(U)}=2n(2n+b), \qquad \nu_n^{(V)}=-2\rho n(2n+b+2). \lab{la_nu_jac} \ea
According to \re{u_expl} this leads to the recurrence coefficients
\be
u_{2n} = 2\,{\frac {n \left( 2\,n-2+a-b \right) }{ \left( a+4\,n \right) 
 \left( a+4\,n-2 \right) }}, \quad u_{2n+1} = {\frac { \left( 2\,n+a \right)  \left( 2\,n+2+b \right) }{ \left( a+4
\,n+2 \right)  \left( a+4\,n \right) }} \lab{u_Jac} \ee
which correspond to those of the generalized Gegenbauer polynomials \cite{Chi}.

For the Laguerre class we have similarly
\ba
&&\lambda_n^{(U)}=-2n + \ve_0, \; \lambda_n^{(V)} = 2\rho n+\ve_1, \nonumber \\
&& \nu_n^{(U)}=2n(2n+b), \: \nu_n^{(V)}=-2\rho n(2n+b+2) \lab{la_nu_lag} \ea
which leads to the recurrence coefficients
\be
u_{2n} = 2n, \quad u_{2n+1} = 2n+2+b. \lab{u_lag} \ee
They correspond to the generalized Laguerre polynomials. For the special case $b=-1$, we obtain the Hermite polynomials.

We have thus described all the admissible classes of symmetric hypergeometric orthogonal polynomials. Note that in contrast to the classification scheme presented in \cite{Zhe_hyp}, only two classes exist for the operator $L$ of the form \re{main_L_x}, namely the $q$-Jacobi and the Jacobi ones. The third class, that of the little -1 Jacobi polynomials, is absent in our classification scheme. The reason for this is the following. The polynomials $U_n(x)$ and $V_n(x)$ introduced in \re{P_UV} are both classical orthogonal polynomials related by the Christoffel transform \re{V_ker}. It is well known that for the ordinary Jacobi and for the little $q$-Jacobi polynomials, the Christoffel transform \re{V_ker} leads to polynomials in the same class but with slightly modified parameters. However, for the little -1 Jacobi polynomials, the Christoffel transform of the polynomials $U_n(x)$ gives non-classical polynomials that do not belong to the little -1 Jacobi family \cite{VZ_little}. This means, in particular, that the compatibility condition \re{add_res} does not hold in this case.

\section{Realizations of the operator $L$}
\setcounter{equation}{0}
The admissible operators $L$ lead to symmetric orthogonal polynomials $P_n(x)$ which are formally bispectral, i.e. they satisfy the ``dual'' eigenvalue equations \re{eig_LP} and \re{rec_sym}. For the classical orthogonal polynomials of the Askey tableau \cite{KLS}, the operator $L$ can be presented in terms of either differential or $q$-difference operators. We now examine this question for the operators $L$ that have arisen in the present study. In contrast to the classical case, here the operator $L$ involves the reflection operator $R$ which means that we deal with the operators of Dunkl type \cite{VZ_little}.

Consider first the case of the generalized Jacobi polynomials. We have the explicit expressions \re{la_nu_jq} for the coefficients $\lambda_n, \nu_n$ of the operator $L$. Introduce the projection operators 
\be
\pi_0 = {1 \over 2} (1+R), \quad    \pi_1 = {1 \over 2} (1-R). \lab{pr_p} \ee
The operator $\pi_0$ projects any function onto the subspace of even functions while the operator $\pi_1$ projects onto the subspace of odd functions.

It is easy to show that the operator $L$ can be presented in the form
\be
L = L_0 \pi_0 + L_1 \pi_1, \lab{L_01} \ee
where the operators $L_0$ and $L_1$ are both 4-difference operators corresponding to the little $q$-Jacobi polynomials. In more details, let us introduce the operators $T_q^+$ and $T_q^-$ defined by the following actions:
\be
T_q^+ f(x) = f(xq), \quad T_q^- f(x) = f(x/q) \lab{T_op} \ee
Then 
\be
L_0 = (-a+bx^{-2}) \left(T_q^+ - \mathcal{I}\right) + (1-x^{-2}) \left(T_q^- - \mathcal{I}\right) + \ve_0 \lab{L0_qj} \ee
and 
\be
L_1 = \rho q (a-bx^{-2}) T_q^+  + \rho q(-1+x^{-2}) T_q^-  + \ve_1 + \rho - a \rho q^2 + \rho(bq^2-1)x^{-2}. \lab{L1_qj} \ee
These expressions for the operators $L_0$ and $L_1$ can easily be verified through their actions on the monomials $x^n$. Indeed, we have
\be
L_0 x^{2n} = \lambda_{2n} x^{2n} + \nu_{2n} x^{2n-2} \lab{L0_xn} \ee
and
\be
L_1 x^{2n+1} = \lambda_{2n+1} x^{2n+1} + \nu_{2n+1} x^{2n-1}. \lab{L1_xn} \ee
Combining, we can present these formulas as
\be
L x^n = (L_0 \pi_0 + L_1 \pi_1)x^n = \lambda_n x^n + \nu_n x^{n-2} \lab{L_tot_monom} \ee
which means that our operator $L$ satisfies the required conditions.

The pure ``classical'' situation arises when $L_0=L_1$. In this case (and only in this case) the reflection operator $R$ vanishes in the  expression of the operator $L$ as seen from the definition of the projection operators \re{pr_p}. It is observed from \re{L0_qj}-\re{L1_qj} that the condition $L_0=L_1$ can be achieved only if 
\be
b=\rho=-q^{-1}, \quad \ve_1 = \ve_0 +(q^{-1}-1)(aq+1). \lab{par_clas_qjac} \ee
Given \re{par_clas_qjac}, the recurrence coefficients \re{u_qj} become
\be
u_n = \frac{q^{n-1}(1-q^n)(1+aq^{n-1})}{(1+aq^{2n})(1+aq^{2n-2})}. \lab{u_big_qjac} \ee
It is possible to identify the orthogonal polynomials with recurrence coefficients given by \re{u_big_qjac} as a special case of the  symmetric big $q$-Jacobi polynomials. Indeed, the big $q$-Jacobi polynomials $P_n(x;\alpha,\beta,\gamma;q)$ depend on 3 parameters $\alpha,\beta,\gamma$ (using the notation of \cite{KLS}) and the recurrence coefficients \re{u_big_qjac} correspond to the choice $a=q\alpha^2, \beta=-\gamma=\alpha$.

Consider now the generalized Gegenbauer polynomials. Applying the same approach, we arrive at the operator \re{L_01} with
\be
L_0= (1-x^2) \partial_x^2 +  \left( (b+1)x^{-1}-(a+1)x  \right) \partial_x + \ve_0 \lab{L0_geg} \ee
and 
\be
L_1= \rho(x^2-1) \partial_x^2 + \rho \left( (a+1)x -(b+1)x^{-1} \right) \partial_x + \ve_1-\rho(a+1) + \rho(b+1) x^{-2} . \lab{L1_geg} \ee
In the special case when $b=-1$ and $\rho=-1, \ve_0=\ve_1=a+1$, we see that  $L_1=L_0$ and hence the operator $L$ becomes classical as it does not contain the reflection operator $R$. This corresponds to the ordinary Gegenbauer polynomials. Indeed, when $b=-1$ one sees from \re{u_Jac} that 
\be
u_n  = \frac{n(n+a-1)}{(2n+a)(2n+a-2)} \lab{gegen_u} \ee 
thus recovering the recurrence coefficients of the  Gegenbauer polynomials $C_n^{(a/2)}(x)$ (in the notation of \cite{KLS}).

For general values of the parameters $a$ and $b$ the operator $L$ is not classical. In \cite{ChGa} the operator $L$ for the generalized Gegenbauer polynomials was presented as a second order Dunkl differential operator.

Consider finally the generalized Laguerre polynomials corresponding to the  choice \re{pr_p}. Again it is easily seen that 
\be
L_0 = \partial_x^2 +\left(-x +(b+1)x^{-1} \right) \partial_x + \ve_0 \lab{L0_lag} \ee
and that
\be
L_1 = -\rho \partial_x^2 + \rho \left(x -(b+1)x^{-1} \right) \partial_x + \ve_1 -\rho +\rho(b+1)x^{-2}. \lab{L1_lag} \ee
When $b=-1$ and $\rho=-1, \ve_0=\ve_1+1$,  the operators $L_1$ and $L_0$ coincide. This corresponds to the ordinary Hermite polynomials. For generic values of $b$ the operator $L$ can be reduced to a second order differential operator of Dunkl type  \cite{ChGa}.

We thus have constructed explicit presentations of the operator $L$ in terms of differential or $q$-difference operators. In the general case, the operator $L$ is non-classical as it contains the reflection operator $R$; only in some special instances (those of the Gegenbauer and Hermite polynomials) can the operator L be reduced to a classical one.

\section{Associated algebras}
\setcounter{equation}{0}
We shall now obtain in our abstract formalism the algebras that underscore the bispectrality of the polynomials that have been identified. In view of the definite parity properties of these polynomials, it will prove pertinent to use the operator of multiplication by the square of the argument $x$:
\be
Y = x^2 \lab{Y_def} \ee
in conjunction with $L$ as generators. Define 
\be
W^{(1)} = L^2 Y + Y L^2 - \Omega LYL, \quad  W^{(2)} = L Y^2 + Y^2 L - \Omega YLY, \lab{W12_def} \ee
where $\Omega=q^2+q^{-2}$ for the little $q$-Jacobi class and $\Omega=2$ for the Jacobi family (including degenerate cases). We shall take $W^{(1)}$ and $W^{(2)}$ to provide the left hand sides of the defining relations of the algebras. 
Consider the following $q$-commutator
\be
[A,B]_q=qAB-q^{-1}BA. \lab{qcom} \ee
Then $W^{(1)}$ and $W^{(2)}$ can be presented  in the form
\be
W^{(1)}= [L,[L,Y]_q]_{q^{-1}}, \quad W^{(2)}= [Y,[Y,L]_q]_{q^{-1}} \lab{W-qcom} \ee
and they obviously become double commutators when $q=1$ and $\Omega=2$. In order to obtain the defining relations in the abstract context, we could determine the action of  $W^{(1)}$ and $W^{(2)}$ on monomials and use the relation \re{main_L_x} to find the identities we are looking for. Alternatively, we can use the explicit representations of the operators $L$ and $Y$ in terms of the differential or $q$-difference operators obtained in the previous section. We shall opt for the latter course in the following.

We start with the case of the generalized little $q$-Jacobi polynomials. There are three independent parameters, $\ve_0, \ve_1, \rho$, in the expression of the coefficients $\lambda_n, \nu_n$ in \re{la_nu_jq}. The orthogonal polynomials $P_n(x)$ do not depend on these constants; only the concrete expression of the operator $L$ will involve these parameters. Their origin can be traced back to formula \re{tL} that gives the operator L as a linear combination (with arbitrary factors) of basic constituents. We can always use the latitude offered by these free constants to simplify the presentations. In the case of the generalized little $q$-Jacobi polynomials it is convenient to put
\be
\rho=-q^{-1}, \; \ve_0=1-a, \; \ve_1= q^{-1}-aq \lab{param_qj} \ee
We then have the relations
\be
[Y,[Y,L]_q]_{q^{-1}}= Y \left(\xi_0   + \xi_1  R \right) \lab{YYL_jj} \ee
and 
\be
[L,[L,Y]_q]_{q^{-1}}= L \left(\xi_0   + \xi_1  R \right)   + \eta Y +\zeta \mathcal{I}, \lab{LLY_jj} \ee
where $R$ is the reflection operator \re{R_def} and where
\ba
&&\xi_0= {1\over 2} (q-1)(b-q^{-1})(q-q^{-1})^2, \quad \xi_1 = {1\over 2} (1-q)(b+q^{-1})(q-q^{-1})^2, \nonumber \\
&&\quad \eta= a (q^2-q^{-2})^2, \quad \zeta= -(a+bq^2)(q^2-q^{-2})^2. \lab{cf_alg_qj} \ea
It is seen that the operators $L$ and $Y$ constitute a special case of a special central extension of the Askey-Wilson algebra $AW(3)$. The central element is the reflection operator $R$ which commutes with both $L$ and $Y$. Koornwinder in \cite{KZ} has introduced and studied the central extension $AW(3,Q_0)$ of the generic Askey-Wilson algebra. Terwilliger has also considered a similar algebraic construction and has called it the ``universal Askey-Wilson algebra'' \cite{T_AW}. In this scheme, the central element is an operator $T$ satisfying the algebraic relation
\be
(T+\alpha \beta )(T+1)=0, \lab{T-alg} \ee
where $\alpha, \beta$ are arbitrary parameters. If one assumes that $\alpha \beta \ne 1$ (the generic situation), by an affine transformation the operator $T$ can be reduced to the reflection operator $R$ which satisfy the algebraic relation $R^2=\mathcal{I}$. We thus see that the algebra \re{YYL_jj} and \re{LLY_jj} can be considered as a special case of $AW(3,Q_0)$ \cite{KZ}, \cite{T_AW}. In the conventional approach, the algebra $AW(3,Q_0)$ is related to the non-symmetric Askey-Wilson polynomials which are Laurent polynomials with nontrivial orthogonality properties and only in the symmetric case do these polynomials reduce to the ordinary Askey-Wilson polynomials (see \cite{KZ} for details). Here, however, the algebra $AW(3,Q_0)$ describes  {\it ordinary} orthogonal polynomials which exhibit a nontrivial bispectral property with respect to the operator $L$.

There is also a relation with a degenerate $q$-Onsager algebra \cite{IT}. Indeed, taking an extra commutators will bring \re{YYL_jj} and \re{LLY_jj} to the relations
\be
[Y, [Y,[Y,L]_q]_{q^{-1}}]= 0, \quad  [L, [L,[L,Y]_q]_{q^{-1}}]= \eta [L,Y] \lab{dg_ons} \ee
which correspond to those of a degenerate $q$-Onsager algebra \cite{IT}. Note that the $q$-Onsager algebra at roots of unity were seen to provide an underpinning for the complementary Bannai-Ito polynomials. \cite{BGV}.

Observe that when $b=-q^{-1}$, the terms with the operator $R$ vanish and our algebra becomes a special case of the ordinary Askey-Wilson algebra. In the previous section we already mentioned that this condition corresponds to a special case of the big $q$-Jacobi polynomials that leads to classical polynomials. In this case, the occurrence of the standard $AW(3)$ algebra is rather natural. In general, it is a central extension of $AW(3)$ that appears.

For the generalized Gegenbauer polynomials we put 
\be
\rho=-1, \; \ve_0= 1-a^2/4, \; \ve_1=-a-a^2/4 \lab{par_gegen} \ee
in \re{L0_lag}-\re{L1_lag}. This brings the following commutation relations
\be
[Y,[Y,L]=-8 Y^2 + 8 Y, \quad [L,[L,Y] =-8 \{L,Y\} + 8 L + 4(b+1) R + \mu \mathcal{I} \lab{gege_alg} \ee 
where $\{L,Y\}=LY+YL$ stands for the anticommutator and where
$$
\mu= 2 a^2 - 4 ab -8a +4b +4.
$$
In this case the algebra \re{gege_alg} is a central extension of the quadratic Jacobi algebra introduced in \cite{GLZ}.

As we already demonstrated, when $b=-1$, the polynomials $P_n(x)$ become the classical Gegenbauer polynomials. With this choice for the parameter $b$, the term containing the reflection operator $R$ disappears in \re{gege_alg} and we return to the ordinary Jacobi algebra.

Finally, for the generalized Hermite polynomials, it is convenient to take the parameters of the operator $L$ as follows:
\be
\ve_0=-1-b/2, \; \ve_1=-2-b/2, \rho=-1 .\lab{par_lag} \ee
We then obtain the relations
\be
[Y,[Y,L]]= 8Y, \quad [L,[L,Y]]= 8L + 4 Y \lab{lie_lag} \ee
which define an algebra isomorphic to $sl_2$. It is interesting to note that even in the case of nonclassical generalized Laguerre polynomials (i.e. when $b \ne -1$), the algebra generated by the operators $Y$ and $L$, is a Lie algebra. This can be achieved by an appropriate choice of the parameters $\ve_0, \ve_1, \rho$. For all other cases of nonclassical polynomials, the appearance of a central extension of the ``classical'' bispectral algebras (like $AW(3)$ and the Jacobi algebra) is inevitable.

\section{Concluding remarks}
\setcounter{equation}{0}
In summary this paper has offered a classification of all positive-definite orthogonal polynomials that are eigenfunctions of an operator $L$ acting as $L x^n = \lambda_n x^n + \nu_n x^{n-2}$ on monomials. This abstract approach has provided a synthetic description of the many symmetric semiclassical orthogonal $q$-polynomials that have been identified recently \cite{GK} (their $q=1$ versions have also be included). For concision, we have taken $0<q \le 1$.

The general form of the non-linear algebras that are associated to the polynomials arising in the classification have also been obtained. Interestingly, these algebras extend those of the Askey-Wilson type in that they involve an involution ($R^2=1$). Algebras with a similar feature have already been encountered in the theory of complementary Bannai-Ito polynomials \cite{GVZ_CBI} and of the dual -1 Hahn polynomials \cite{TVZ_Hahn}. It has recently been found that in these cases, the corresponding algebras can be embedded in $q$-Onsager algebras for $q$ a root of unity \cite{BGV}. It would be interesting to provide a similar connection between tridiagonal pairs and the algebras identified in this paper. We plan on returning to this question in the future.

\section*{Acknowledgments}

\noindent 
S.~T., L.~V. and A.~Z. have appreciated the hospitality of the School of Mathematical Sciences of the Shanghai University where part of the research presented in this paper has been realized. Luc Vinet holds a Senior Chair Visiting Professorship at the Shanghai Jiao Tong University which he gratefully acknowledges. The research of ST is supported by  JSPS KAKENHI (Grant Numbers 25400110), that of LV by NSERC and that of GY by NNSF (grant \# 11371251).


\bb{99}

\bi{Al-Salam} W.~A.~Al-Salam, 
{\it Characterization theorems for orthogonal polynomials}, 
in: P. Nevai (ed.), {\it Orthogonal Polynomials: Theory and Practice}, NATO Adv. Sci. Inst. Ser. C: Math. Phys. Sci., vol. 294, Kluwer Academic Publ., Dordrecht, 1990, pp. 1--24.  

\bi{AW} R.~Askey and J.~Wilson, 
{\it Some basic hypergeometric orthogonal polynomials that generalize Jacobi polynomials}, 
Mem. Amer. Math. Soc. {\bf 319} (1985) 1--55.

\bi{BGV} P.~Baseilhac, A.~M.~Gainutdinov and T.~T.~Vu,
{\it Cyclic tridiagonal pairs, higher order Onsager algebras and orthogonal polynomials}, 
arXiv:1607.00606.

\bi{ChGa} Y.~Ben~Cheikh and M.~Gaied,
{\it Characterization of the Dunkl-classical symmetric orthogonal polynomials}, 
Appl. Math. Comput. {\bf 187} (2007) 105--114.

\bi{ChGaZa} Y.~Ben~Cheikh, M.~Gaied and A.~Zaghouani,
{\it $q$-Dunkl-classical $q$-Hermite type polynomials}, 
Georgian Math. J. {\bf 21} (2014) 125--137.

\bi{Chi} T.~Chihara, 
{\it An Introduction to Orthogonal Polynomials}, 
Gordon and Breach, NY, 1978.

\bi{Duran} A.~J.~Dur\'an, 
{\it Orthogonal polynomials satisfying higher-order difference equations}, 
Constr. Approx. {\bf 36} (2012) 459--486.

\bi{GVZ_CBI} V.~X.~Genest, L.~Vinet and A.~Zhedanov, 
{\it Bispectrality of the complementary Bannai-Ito polynomials}, 
SIGMA {\bf 9} (2013) 018, 20 pages; arXiv:1211.2461.

\bi{GK} A.~Ghressi and L.~Kh\'eriji,
{\it The symmetrical $H_q$-semiclassical orthogonal polynomials of class one}, 
SIGMA {\bf 5} (2009), 076, 22 pages.

\bi{GLZ} Y.~A.~Granovskii, I.~Lutzenko and A.~Zhedanov,
{\it Mutual integrability, quadratic algebras and dynamical symmetry},
Ann. Physics {\bf 217} (1992) 1--20.

\bi{Grun} F.~A.~Gr\"unbaum and L.~Haine, 
{\it Bispectral Darboux transformations: an extension of the Krall polynomials}, 
Int. Math. Res. Not. (1997) 359--392.

\bi{Ismail} M.~E.~H.~Ismail, 
{\it Classical and Quantum Orthogonal Polynomials in One Variable},
Encyclopedia of Mathematics and its Applications {\bf 98}, Cambridge University Press, Cambridge, 2005.

\bi{IT} T.~Ito and P.~Terwilliger, 
{\it The augmented tridiagonal algebra}, 
Kyushu J. Math. {\bf 64} (2010) 81--144; arXiv:0904.2889.

\bi{Khol} A.~N.~Kholodov, 
{\it The umbral calculus and orthogonal polynomials}, 
Acta Appl. Math. {\bf 19} (1990) 1--54. 

\bi{KLS} R.~Koekoek, P.~A.~Lesky and R.~F.~Swarttouw,
{\it Hypergeometric Orthogonal Polynomials and Their $q$-Analogues}, 
Springer Monographs in Mathematics, Springer-Verlag, Berlin, 2010.

\bi{KZ} T.~Koornwinder, 
{\it The relationship between Zhedanov algebra $AW(3)$ and the double affine Hecke algebra in the rank one case}, SIGMA {\bf 3} (2007) 063, 15 pages.

\bi{NSU} A.~F.~Nikiforov, S.~K.~Suslov and V.B. Uvarov, 
{\it Classical Orthogonal Polynomials of a Discrete Variable},
Springer-Verlag, Berlin, 1991.

\bi{Roman} S.~Roman, 
{\it The theory of the umbral calculus. I }, 
J. Math. Anal. Appl. {\bf 87} (1982) 58--115. 

\bi{SVZ} V.~Spiridonov, L.~Vinet and A.~Zhedanov, 
{\it Bispectrality and Darboux transformation in the theory of orthogonal polynomials},
in:  J.~Harnad and A.~Kasman (eds.), {\it The Bispectral Problem}, CRM Proc. Lecture Notes, vol. 14, Amer. Math. Soc., Providence, 1998, pp. 111--122.

\bi{SZ_root} V.~Spiridonov and A.~Zhedanov, 
{\it Zeros and orthogonality of the Askey-Wilson polynomials for $q$ a root of unity},
Duke Math. J. {\bf 89} (1997) 283--305.

\bi{Ter} P.~Terwilliger, 
{\it Two linear transformations each tridiagonal with respect to an eigenbasis of the other},
Linear Algebra Appl. {\bf 330} (2001) 149--203.

\bi{T_AW} P.~Terwilliger, 
{\it The universal Askey-Wilson algebra}, 
SIGMA {\bf 7} (2011) 069, 24 pages.

\bi{TVZ_Hahn} S.~Tsujimoto, L.~Vinet and A.~Zhedanov, 
{\it Dual $-1$ Hahn polynomials: ``classical'' polynomials beyond the Leonard duality}, 
Proc. Amer. Math. Soc. {\bf 141} (2013) 959--970; arXiv:1108.0132.

\bi{VYZ} L.~Vinet, O.~Yermolayeva and A.~Zhedanov, 
{\it A method to study the Krall and $q$-Krall polynomials}, 
J. Comp. Appl. Math. {\bf 133} (2001) 647--656.

\bi{VZ_Bochner} L.~Vinet and A.~Zhedanov, 
{\it Generalized Bochner theorem: characterization of the Askey-Wilson polynomials},
J. Comp. Appl. Math. {\bf 211} (2008) 45--56.

\bi{VZ_little} L.~Vinet and A.~Zhedanov, 
{\it A `missing' family of classical orthogonal polynomials}, 
J. Phys. A: Math. Theor. {\bf 44} (2011) 085201, 16 pages; arXiv:1011.1669.

\bi{VZ_Newton} L.~Vinet and A.~Zhedanov, 
{\it Hypergeometric orthogonal polynomials with respect to Newtonian bases}, 
SIGMA {\bf 12} (2016) 048, 14 pages.

\bi{Zhe_hyp} A.~Zhedanov, 
{\it Abstract ``hypergeometric'' orthogonal polynomials}, 
arXiv:1401.6754.

\bi{ZheAW}  A.~Zhedanov, 
{\it Hidden symmetry of Askey-Wilson polynomials}, 
Theoret. Math. Phys. {\bf 89} (1991) 1146--1157.

\bi{Z_umbral} A.~Zhedanov, 
{\it Umbral ``classical'' polynomials}, 
J. Math. Anal. Appl. {\bf 420} (2014) 1354--1375.

\eb

\end{document}